\documentclass[11pt]{article}
\advance \topmargin by -\headheight    \advance \topmargin by -\headsep
\usepackage{times}
\usepackage{xcolor}
\usepackage{soul}
\usepackage[utf8]{inputenc}
\usepackage[small]{caption}

\usepackage{latexsym} 

\newcommand{\pile}[2][c]{\begin{tabular}{@{}#1@{}}#2\end{tabular}}

\newcommand{\apile}[2][c]{\begin{array}{@{}#1@{}}#2\end{array}}
\newcommand{\frc}[2]{{}^{#1\!\!}/_{\!#2}}

\newcommand{\zz}[2][l]{\makebox[0pt][#1]{\raisebox{0pt}[0pt][0pt]{#2}}}
\newcommand{\pmat}[1]{\left(\begin{matrix}#1\end{matrix}\right)}

\newcommand{\bb}{\mathbf}
\newcommand{\DS}{\displaystyle}
\newcommand{\mc}{\multicolumn}
\newcommand{\itize}[1]{\par\begin{itemize}\itemsep=0pt\topsep=-0ex\parskip=0pt#1\end{itemize}}

\newcommand{\boldparagraph}[1]{\par\vspace{.25ex}\noindent{\bf #1.}}

\usepackage{amsmath}
\usepackage{amssymb}
\usepackage{mathtools}
\usepackage{relsize}
\usepackage{bm}
\usepackage{color}
\usepackage{tikz}
\usepackage{enumitem}
\def\myeqno#1{\refstepcounter{equation}\label{#1}\eqno{(\theequation)}}
\newcounter{thm} %
\newcommand{\nextthmno}{\refstepcounter{thm}{\bf \arabic{thm}}}
\newcommand{\newthm}[1][Theorem]{\par\noindent{\bf #1} \nextthmno{\bf .}}
\newcommand{\newthmx}[1][Theorem]{\par\noindent{\bf #1} \nextthmno}
\newcommand{\Proof}[1][Proof]{\par\noindent{\bf #1. }}
\newcommand{\QED}{\mbox{~}\hfill\rule{1ex}{2ex}}%

\newcommand{\Prob}{{\ensuremath{\mathrm{Pr}}}}
\newcommand{\PHT}{\mathbf{PHT}}

\emergencystretch=\maxdimen
\title{Random Walk Laplacian and Network Centrality Measures\thanks{This research was supported in part by NSF grants
1319749 \& 1460620.%
}}

\iftrue
\author{
 \pile{Daniel Boley\\Univ. of Minnesota}
 ~~
 \pile{Alejandro Buendia\\Columbia Univ.}
 ~~
 \pile{Golshan Golnari\\3M}
}
\fi

\begin{document}
\newsavebox{\GGaug}
\sbox{\GGaug}{
\begin{tikzpicture}[scale=0.13]
\tikzstyle{every node}+=[inner sep=0pt]
\draw [black] (17.2,-23.3) circle (3);
\draw (17.2,-23.3) node {$v_1$};
\draw [black] (38.1,-7.3) circle (3);
\draw (38.1,-7.3) node {$v_2$};
\draw [black] (57.1,-20.4) circle (3);
\draw (57.1,-20.4) node {$v_3$};
\draw [black] (52.8,-44) circle (3);
\draw (52.8,-44) node {$v_4$};
\draw [black] (22.9,-44) circle (3);
\draw (22.9,-44) node {$v_5$};
{\color{black}
\draw [black] (70,-32.8) circle (3);
\draw (70,-32.8) node {$v_6$};
}
\draw [black] (35.92,-9.36) arc (-47.38425:-57.74401:112.722);
\fill [black] (35.92,-9.36) -- (34.99,-9.53) -- (35.67,-10.27);
\draw (31.00,-15,00) node [below] {$.340$}; %
\draw [black] (19.742,-21.711) arc (118.19724:1.45013:22.539);
\fill [black] (52.92,-41) -- (53.4,-40.19) -- (52.4,-40.22);
\draw (40.00,-21.00) node [above] {$.510$}; %
\draw [black] (49.978,-45.016) arc (-72.35704:-107.64296:40.015);
\fill [black] (25.72,-45.02) -- (26.33,-45.73) -- (26.64,-44.78);
\draw (37.85,-47.4) node [below] {$.425$}; %
\draw [black] (40.57,-9) -- (54.63,-18.7);
\fill [black] (54.63,-18.7) -- (54.26,-17.83) -- (53.69,-18.65);
\draw (49.85,-13.35) node [above] {$.425$}; %
\draw [black] (49.95,-43.063) arc (-109.12986:-131.22277:92.147);
\fill [black] (19.42,-25.31) -- (19.7,-26.22) -- (20.36,-25.46);
\draw (34.50,-34.00) node [below] {$.425$}; %
\draw [black] (24.242,-41.317) arc (151.30281:97.91308:40.398);
\fill [black] (24.24,-41.32) -- (25.06,-40.86) -- (24.19,-40.38);
\draw (37.00,-26.00) node [above] {$.850$}; %
\draw [black] (40.215,-9.427) arc (43.02289:0.63375:47.04);
\fill [black] (52.86,-41) -- (53.35,-40.2) -- (52.35,-40.21);
\draw (39.00,-14.00) node [right] {$.255$}; %
\draw [black] (22.1,-41.11) -- (18,-26.19);
\fill [black] (18,-26.19) -- (17.73,-27.1) -- (18.69,-26.83);
\draw (19.28,-34.18) node [left] {$.170$}; %
\draw [black] (18.305,-20.513) arc (154.65789:100.21385:23.149);
\fill [black] (18.31,-20.51) -- (19.1,-20) -- (18.2,-19.58);
\draw (23.4,-11.54) node [above] {$.170$}; %
\draw [black] (55.299,-22.799) arc (-38.29642:-72.48769:60.999);
\fill [black] (55.3,-22.8) -- (54.41,-23.12) -- (55.2,-23.74);
\draw (43.82,-35.7) node [below] {$.680$}; %
\color{black}
\draw [black] (67,-32.753) arc (-91.44947:-108.95007:156.863);
\fill [black] (67,-32.75) -- (66.21,-32.23) -- (66.19,-33.23);
\draw (40.00,-30.00) node [below] {$.150$}; %
\draw [black] (41.09,-7.535) arc (82.98946:19.73474:34.206);
\fill [black] (69.11,-29.94) -- (69.31,-29.01) -- (68.37,-29.35);
\draw (59.39,-14.27) node [above] {$.150$}; %
\draw [black] (59.26,-22.48) -- (67.84,-30.72);
\fill [black] (67.84,-30.72) -- (67.61,-29.81) -- (66.91,-30.53);
\draw (62.42,-27.08) node [below] {$.150$}; %
\draw [black] (55.31,-42.36) -- (67.49,-34.44);
\fill [black] (67.49,-34.44) -- (66.54,-34.45) -- (67.09,-35.29);
\draw (62.51,-38.9) node [below] {$.150$}; %
\draw [black] (67.12,-33.641) arc (-73.91148:-79.33642:448.199);
\fill [black] (67.12,-33.64) -- (66.21,-33.38) -- (66.49,-34.34);
\draw (37.00,-42.00) node [below] {$.150$}; %
\draw (72.5,-35.5) node [below] %
  {{\pile[ll]{\scalebox{.75}[1]{\large$\swarrow\!\!\downarrow\!\!\searrow$}\,0.200 \\\hspace*{-1ex}\sl \fbox{$v_{1\cdots 5}$}}}};
\end{tikzpicture}%
}

\maketitle

\centerline{Abstract}

Random walks over directed graphs are used to model activities in many domains, such as social networks, influence propagation, and Bayesian graphical models. They are often used to compute the importance or centrality of individual nodes according to a variety of different criteria.  Here we show how the pseudoinverse of the ``random walk" Laplacian can be used to quickly compute measures such as the average number of visits to a given node and various centrality and betweenness measures for individual nodes, both for the network in general and in the case a subset of nodes is to be avoided.  We show that with a single matrix inversion it is possible to rapidly compute many such quantities.
\medskip

\section{Introduction}
\boldparagraph{Background}
We consider random walks over a strongly connected directed graph, 
corresponding to a recurring Markov chain.   A graph with $n$ vertices
can be fully described by its $n \times n$ adjacency matrix $A$ whose
$ij$-th entry $a_{ij}$ is the weight on the directed edge from $i$ to $j$,
or zero if there is no such edge.
Let $\bb 1$ be the vector of all ones.
We let $\bm d = A \cdot \bb 1$ be the vector of out-degrees, where $d_i = \sum_j a_{ij}$, and
$D = \textsc{Diag}\bm(d)$ be the diagonal matrix with the entries of
$\bm d$ on the diagonal.
Then the corresponding Markov chain has transition probability matrix
$P = D^{-1} A$.
It is well known that if the graph is strongly connected, then 
$P$ has a simple eigenvalue $\lambda = 1$ with corresponding
eigenvector with all positive entries $\bm \pi = (\pi_1,\ldots,\pi_n)^T$,
which can be scaled to unit length in the 1-norm, i.e. the vector of stationary probabilities for the
random walk.
We let $\bm\Pi = \textsc{Diag}(\bm\pi)$ be the corresponding diagonal matrix.
There are several Laplacians that can be formed from this matrix
as noted in \cite{Boley11}:
$${
\apile[l@{\,\,=\,\,}l@{\quad}l]{
L & I - P & \mbox{normalized Laplacian} \\
\bb L & \bm\Pi (I - P) & \mbox{random walk Laplacian} \\
\cal L & D-A & \mbox{combinatorial Laplacian}
}
}\myeqno{laplacians}$$
It is well known that these Laplacians play critical roles with respect to undirected
graphs.  For undirected graphs, we have $\bm d = \bm \pi \cdot (\mbox{volume})$,
where $\mbox{volume}$ is twice the number of edges in the graph, and hence the 
combinatorial Laplacian is just a scalar multiple of the random walk Laplacian.
For a directed graph, these ``Laplacians'' are no longer symmetric, but
\cite{Boley11} showed how they can still be used to compute interesting
properties of the graph.

Golnari et al.\ (see e.g. \cite{Golnari15}) introduced the third-order
{\it random walk fundamental tensor} $\bm{N} = \{\bm{N}(i,j,k)\}_{i,j,k=1}^n$ where each entry
$\bm N(i,j,k)$ is the expected number of visits to intermediate node $j$
on all walks starting from node $i$ before reaching target node $k$.
This is an extension of the so-called 
fundamental matrix $N$ of an absorbing Markov chain \cite{Grinstead06},
which
is a matrix whose $ij$-th entry gives the expected number of random walk passages through node $j$ when starting a
random walk from node $i$ before being absorbed.
Each slice $\bm N (:,:,k)$ of the tensor is the fundamental
matrix for the modified Markov chain where node $k$ is turned into an absorbing node.
This tensor can be used to quickly compute several centrality measures such as random walk closeness \(\sum_{i} H_i^{[k]} = \sum_{i,j}\bm{N}(i,j,k)\) for a node $k$ \cite{Noh04}, and random walk betweenness \(\sum_{i\not=j,k\not=j}\Prob({i  \rightarrow   j   \rightarrow   k})\)
(i.e., the sum of probabilities of a random walk passing an intermediate node $j$ averaged over all starting nodes $i$ and ending nodes $k$) \cite{Newman05,Kang11}. %

\boldparagraph{Outline} %
In this paper we review the random walk Laplacian and show that we can compute its pseudoinverse
using one matrix inversion (\(\mathcal{O}(n^3)\) complexity) plus other computations of lower complexity, 
many of which can be found in \cite{Boley11}. We then show how this pseudoinverse can be used to efficiently
calculate many useful properties of the graph such as the fundamental tensor \cite{GolnariXX},
which leads to convenient calculation of hitting times, absorption probabilities, and centrality measures, in
addition to similar measures conditioned on avoiding certain nodes.

After the preliminaries, we sketch some properties of the random walk Laplacian, formally define the random walk fundamental tensor, and show how the tensor can be used to compute properties of the underlying graph based on its correspondence to the random walk.

\boldparagraph{Related Work}
This work was motivated by a long literature of centrality and betweenness measures \cite{Brandes01,Barthelemy04,Newman05,Fouss07,Kang11,Mavro15}, which have
traditionally been difficult to compute.  Random walk-based centrality measures have improved on other notions of centrality by accounting for propagation through all possible paths between a source and target. This is in contrast to previous metrics that have ranked agents exclusively by geodesics or by maximum flow through a particular choice of idealized paths \cite{Stephenson89,Freeman91}. Random walk probabilities of passages through individual nodes are also
used in influence propagation, such as in trust mechanisms \cite{Hopcroft07,Hang10,Liu16}.

\section{Random Walk Laplacian and its Pseudoinverse} %

We consider a strongly connected directed graph
\(\mathcal{G} = (\mathcal{V}, \mathcal{E}, A)\),
where $|\mathcal{V}|= n$.
This digraph
can be modeled by an irreducible Markov chain with stationary probability vector \(\bm{\pi}\).
By strongly connected we mean that there is a path from any node to any other
node within the directed graph.
For simplicity we assume the graph is
unweighted, though the development can easily be extended to the weighted case.
The probability transition matrix for the Markov chain is calculated as \(P = D^{-1}A\), in which
$A$ is the adjacency matrix and \(D=\textsc{Diag}(A\cdot\bb 1)\) is the diagonal matrix of vertex out-degrees.

The probability transition matrix for the random walk, treated as a 
Markov chain, is
$${
P = 
\begin{pmatrix}
P_{\alpha,\alpha} & \mathbf{r} \\
\mathbf{s}^T & \tau
\end{pmatrix} 
=
\begin{pmatrix}
P_{\alpha,\alpha} & P_{\alpha,{}n} \\
P_{{}n,\alpha} & P_{{}n{}n}
\end{pmatrix}  ,
}\myeqno{equ:P}$$
where $\alpha = \{1,\ldots,n-1\}$ is the index set corresponding to all but the last vertex (in their arbitrary order),
and $\delta=\{n\}$ is a chosen absorbing state.
This matrix is row-stochastic (i.e., $P\cdot \bb 1 = \bb 1$), and the
vector of recurring probabilities $\bm\pi = (\pi_1;\ldots;\pi_n)$ satisfies
$\bm\pi^T P = \bm\pi^T $, $\pi_1>0,\ldots,\pi_n>0$, $\|\bm \pi\|_1=1$.
Here we use ``;'' to denote vertical concatenation {\it \`a la} Matlab.
If we were to assume no self-loops, $\tau$ would be zero.

The matrix $P$ has a simple eigenvalue equal to 1 with right eigenvector
$\bb 1$ and left eigenvector $\bm\pi^T$:
$${ P\cdot \bb 1 = \bb 1 \quad \quad \bm\pi^T P = \bm\pi^T}$$
The left eigenvector is scaled by $\|\bm\pi\|_1 = \sum_j \pi_j = 1$ so that 
it is the vector of stationary probabilities for the recurring Markov chain.
The fact that this eigenvalue is simple is implied by the assumption that the
underlying graph is strongly connected.

The normalized Laplacian is $L = I-P$ and the random walk Laplacian is
$\bb L = \bm\Pi(I-P)$.   Both of these have rank $n-1$ and nullity 1.
The left and right annihilating vectors for $L$ are exactly the left and right eigenvectors
of $P$ corresponding to eigenvalue 1, namely $\bb 1$ and $\bm \pi$, respectively.
The annihilating vectors for $\bb L $ are both $\bb 1$ (on both sides):
$${ \bb L\cdot \bb 1 = \bm 0 \quad \quad \bm 1^T \bb L = \bm 0^T .}$$
This implies the Moore-Penrose inverse of $\bb L $ is the same as its
Drazin inverse, making it particularly easy to write down its pseudoinverse.
These properties are summarized in \cite{Boley11}.

We use the following lemma based on \cite[Lemma 1]{Boley11} regarding the
Moore-Penrose inverse of matrices with nullity 1.

\newthm[Lemma]\label{inv}
Suppose we have an invertible $(n-1)\times (n-1)$ matrix $L_{\alpha,\alpha}$ and
two $(n-1)$-vectors $\bb u, \bb v$.  Then there is a unique
$n \times n$ matrix $L$ 
with nullity 1 (rank $n-1$)
of the form
$${
\apile[lll]{
L &=& \pmat{L_{\alpha,\alpha} & \bm l_{\alpha,{}n}\cr \bm l_{{}n,\alpha}^T   & l_{nn}}
\\ &=&
\pmat{L_{\alpha,\alpha} & - L_{\alpha,\alpha} \bb v
\cr - \bb u^T L_{\alpha,\alpha} & \bb u^T L_{\alpha,\alpha} \bb v}
\\ &=&
\pmat{I_{n-1} \cr -\bb u^T} L_{\alpha,\alpha} \pmat{ I_{n-1} ,\, -\bb v}
}}\myeqno{fillL}$$
having left and right annihilating vectors satisfying
$${(\bb u^T   , 1) L = 0, \quad  L (\bb v; 1) = 0}.\myeqno{annihil}$$ %

Furthermore the Moore-Penrose pseudoinverse $M = L^+$ can be written
in terms of $L_{\alpha,\alpha},\, \bb u,\, \bb v$ as follows:
$${
\apile[lll]{
M &=& \pmat{M_{\alpha,\alpha} & \bm m_{\alpha,{}n}\cr \bm m_{{}n,\alpha}^T   & m_{nn}}
\\ &=&
\pmat{R_{\bb v} \cr \frac{-1}{1+\bb v^T\bb v}\bb v^T} L_{\alpha,\alpha}^{-1} 
  \pmat{R_{\bb u} ,\, \frac{-1}{1+\bb u^T\bb u} \bb u}
}}\myeqno{partitionM}$$
(where $M_{\alpha,\alpha}$ is the upper-left $(n-1)\times(n-1)$ block of $M$)
with
$${
\apile[r@{\,=\,}l@{\quad\quad\quad}r@{\,=\,}l]{
M_{\alpha,\alpha}  & R_{\bb v} 
                L_{\alpha,\alpha}^{-1}
                R_{\bb u}
& \bm m_{\alpha,{}n} & - M_{\alpha,\alpha} \bb u
\\
\bm m_{{}n,\alpha}^T & - \bb v^T M_{\alpha,\alpha} & m_{nn} & \bb v^T M_{\alpha,\alpha} \bb u, 
}
}\myeqno{fillM}$$
where $R_{\bb u} = (I_{n-1} + \bb u \bb u^T )^{-1} = (I_{n-1} - \frac{1}{1+\bb u^T\bb u} \bb u \bb u^T   )$,
and $R_{\bb v} = (I_{n-1} + \bb v \bb v^T )^{-1} = (I_{n-1} - \frac{1}{1+\bb v^T\bb v} \bb v \bb v^T   )$.

\Proof[Proof (Sketch)]

Equation (\ref{fillL}) is a simple consequence of (\ref{annihil}). 
Regarding $M$, it can be
verified by direct calculation and some simplification that 
$ML$ and $LM$ are symmetric, $LML=L$, and $MLM=M$, satisfying
the 
conditions for the
Moore-Penrose pseudoinverse.
\QED

\newthm[Corollary] \label{thm:L11} \cite{Boley11}
Suppose $M$ is an $n \times n$ matrix partitioned as in (\ref{fillM}) with
the upper $(n-1)\times (n-1)$ block $M_{\alpha,\alpha}$ invertible, and satisfying
$(\bb v^T   , 1) M = 0$ and $M (\bb u; 1) = 0 $.
Then the upper-left
$(n-1)\times(n-1)$ block of $L = M^+$ satisfies
\begin{equation} \label{equ:L11}
\apile[lll]{
L_{\alpha,\alpha}^{-1} &=& R_{\bb v}^{-1} M_{\alpha,\alpha} R_{\bb u}^{-1}
\\ &=&
(I_{n-1} +\bb v \bb v^T   ) M_{\alpha,\alpha} (I_{n-1} +\bb u \bb u^T   )
\\ &=&
(I_{n-1}  \,,\, -\bb v )
\pmat{M_{\alpha,\alpha} & \bm m_{\alpha,{}n}\cr \bm m_{{}n,\alpha}^T   & m_{nn}} 
\pmat{I_{n-1}  \cr - \bb u^T   } .
}\end{equation}
\QED

This corollary leads directly to an algorithm to compute the random walk Laplacian and its pseudoinverse
using only one matrix inverse.

\rule[-1in]{0pt}{1in}\vspace{-1in}%
\newthmx[Algorithm]{\bf{ Compute Moore-Penrose pseudoinverse of random walk Laplacian.}}\label{algM}
\begin{enumerate}
\item Compute normalized Laplacian $L=I-P$ partitioned as in (\ref{fillL}).
\item Compute some representation of the inverse of the upper $(n-1)\times (n-1)$ part
$L_{\alpha,\alpha}$\footnote{The inverse could be represented in terms of its LU factors, but
later we need the entries of the inverse itself.}.
\item Solve the linear system
$(\pi_1,\ldots,\pi_{n-1}) $ $=$ $ - L_{\alpha,\alpha}^{-1}\bm l_{\alpha,{}n} \pi_n$,
where $\pi_n$ is scaled after the fact to satisfy $\|\bm \pi\|_1=1$.
\item Form diagonal matrix $\bm \Pi = \textsc{Diag}(\bm \pi)$ and
random walk Laplacian $\bb L = \bm \Pi \cdot L = \bm \Pi (I-P)$, itself partitioned
as in (\ref{fillL}).
\item Compute the inverse of the upper-left block of $\bb L$,
namely $\bb L_{\alpha,\alpha}^{-1} = (I - P_{\alpha,\alpha})^{-1} \bm \Pi_1^{-1}$
using the previously computed inverse, where
$\bm \Pi_1 = \textsc{Diag}(\pi_1,\ldots,\pi_{n-1})$.
\itize{\item Elementwise: $[\bb L^{-1}]_{ij} = [L^{-1}]_{ij} / \pi_j$.}
\item Compute desired pseudoinverse $\bb M$ according to Lemma \ref{inv}: %
$${
\apile[lll]{
\bb M &=& \pmat{\bb M_{\alpha,\alpha} & \bb m_{\alpha,{}n}\cr \bb m_{{}n,\alpha}^T   & \mathrm m_{nn}}
\\ &=&
\pmat{R_{\bb 1} \cr \frac{-1}{n}\bb 1^T} \bb L_{\alpha,\alpha}^{-1} 
  \pmat{R_{\bb 1} ,\, \frac{-1}{n} \bb 1},
}}\myeqno{fillBM}$$
where we have used the fact that the left and right annihilating vectors for
$\bb L$ and $\bb M$ are both $\bb 1_n$, so the vectors in (\ref{annihil}) are $\bb u = \bb v = \bb 1_{n-1}$. We also have
the identity $R_{\bb 1}\bb 1 = (I_{n-1} - \frac{1}{n} \bb 1 \bb 1^T ) \bb 1 = \frac{1}{n} \bb 1$.
\\This can be computed step-by-step as follows:
\itize{
\item[a.] Compute $\bb b = \frc{1}{n} \bb L_{\alpha,\alpha}^{-1} \bb 1$ and $\bb c^T = \frc{1}{n}\bb 1^T \bb L_{\alpha,\alpha}^{-1}.$
\item[b.] Compute $\bb M_{\alpha,\alpha} = \bb L_{\alpha,\alpha}^{-1} - \bb b \cdot \bb 1^T - \bb 1 \cdot \bb c^T + \frac{\bb c^T \bb 1}{n} \bb 1
\bb 1^T$,\\ or elementwise
$[\bb M_{\alpha,\alpha}]_{ij} = [\bb L_{\alpha,\alpha}^{-1}]_{ij} - [\bb b]_i - [\bb c]_j + \frac{\bb c^T \bb 1}{n}$.
\item[c.] Compute $\bb m_{\alpha,{}n} = \bb b - \frac{\bb 1^T \bb b}{n} \bb 1$,
$\bb m_{{}n,\alpha}^T = \bb c^T - \frac{\bb c^T \bb 1}{n} \bb 1^T$,
$\mathrm m_{nn}= \frac{\bb c^T \bb 1}{n} = \frac{\bb 1 ^T \bb b}{n}$.
}
\end{enumerate}
We remark that in the case of undirected graphs, $\bm \pi$ is a multiple of the vector of vertex degrees, hence step 2 would be free, but we would still need the inverse in step 5.

\section{Fundamental Tensor}
\boldparagraph{Fundamental Matrix}
If we set $\mathbf{s}^T=0, \tau=1$ in (\ref{equ:P}), we turn the recurring Markov chain
into an absorbing Markov chain 
with corresponding probability transition matrix $\widetilde P$,
where the last node $k=n$ is a single absorbing node.
In this case the {\it fundamental matrix} $N$ \cite{Grinstead06} is the matrix
$${ 
\bm N(\alpha,\alpha,n) = L_{\alpha,\alpha}^{-1}  = (I - P_{\alpha,\alpha})^{-1},
}\myeqno{Nlap}$$
whose $ij$-th entry is the average or expected 
number of passages through node $j$ in random walks starting from node $i$ which are absorbed by
node $k$.

\boldparagraph{Fundamental Tensor}
We would like to compute the {\it fundamental tensor} $\bm N = \{\bm N(i,j,k)\}_{i,j,k=1,\ldots,n}$
where $\bm N(i,j,k)$ is the average number of passages through
intermediate node $j$ for random walks starting from $i$ before reaching $k$.
This would be the fundamental matrix in a Markov chain in which node $k$ is made
absorbing instead of node $n$.
In Algorithm \ref{algM} we have shown how to compute the pseudoinverse of the random walk
Laplacian from the fundamental matrix $N = \bm N(\alpha,\alpha,n)$ obtained when node $n$ is made absorbing.
Here $\alpha=\{1,\ldots,n-1\}$ denotes the index set selecting all the nodes
except the $n$-th.
We would now like to do the reverse: obtain the fundamental matrix $N$ from the pseudoinverse
$\bb M$ of the random walk
Laplacian without any more matrix inversions.  Since the choice of absorbing node $n$ is
arbitrary, this method will also yield a method to compute the fundamental
matrix $\bm N(:,:,k)$ obtained when any other node $k \not = n$ is made absorbing.  In this way
we can fill in the entire fundamental tensor.  

\boldparagraph{Generalization of Fundamental Tensor}
In later sections below, the submatrix of $N=\bm N(\alpha,\alpha,n)$ consisting of rows indexed by
index set $\beta$ and columns by index set $\gamma$ is denoted by $N_{\beta,\gamma} = \bm N(\beta,\gamma,n)$.
We also later discuss other tensors of visit counts such as
$${\apile[l@{~}l@{~}l@{~}l]{
\mbox{(a)} &  \bm N (i,j,\{\gamma,n\}) &=& \mbox{\parbox[t]{4in}{expected number of visits to $j$ before reaching any node in $\{\gamma,n\}$, starting from $i$}}
\\
\mbox{(b)} & \bm N(i,j,n,\gamma) &=& \mbox{\parbox[t]{4in}{expected number of visits to $j$ before reaching $n$ starting from $i$, counting only walks that avoid $\gamma$}}
}
}$$ 

Corollary \ref{thm:L11} provides the formula to go from $\bb M$ to $N=\bm N(\alpha,\alpha,n)$.
From (\ref{fillM}), it follows that
$${
N = (I-P_{\alpha,\alpha})^{-1} = L_{\alpha,\alpha}^{-1} = \bb L_{\alpha,\alpha}^{-1} \bm \Pi_1
= R_{\bb 1}^{-1} \bb M_{\alpha,\alpha} R_{\bb 1}^{-1} \bm \Pi_1
}$$
A little algebra yields
$${
\apile[lll]{
N &=& (I+\bb 1 \bb1^T) \bb M_{\alpha,\alpha} (I+\bb 1 \bb1^T)
\\&=&
\pmat{I ,\, -\bb 1} \pmat{\bb M_{\alpha,\alpha} & \bb m_{\alpha,{}n} \cr \bb m_{{}n,\alpha}^T & \bm m_{nn}} \pmat{I \cr -\bb 1^T} \bm \Pi_1
\\&=&
\pmat{I ,\, -\bb 1} \bb M \pmat{I \cr -\bb 1^T} \bm \Pi_1,
}}$$
or elementwise
$N(i,j) = (\bm m_{ij} - \bm m_{nj} - \bm m_{in} + \bm m_{nn}) \pi_j$ for $1 \le i,j \le n-1$.
This yields the values for the $ij$-th entries in the $n$-th slice of the tensor, $\bm N(i,j,n)$,
including the values when $i=n$ or $j=n$.
We are free to reorder the nodes because the left and right annihilating vectors
of $\bb L$ remain the same regardless of the order (namely, the vector of all ones).
By reordering the nodes, we obtain the corresponding formula for every other slice
in terms of the entries of $\bb M$:
$${
\bm N(i,j,k) =  (\bm m_{ij} - \bm m_{kj} - \bm m_{ik} + \bm m_{kk}) \pi_j
}\myeqno{equ:slice}$$
We can compute the entire fundamental tensor $\bm N$ using Algorithm \ref{algM}
plus formula (\ref{equ:slice}) in $O(n^3)$ time, requiring only one matrix inverse and
averaging out to constant time per entry $\bm N(i,j,k)$.

To obtain the generalizations $\bm N(\beta,\beta,n,\gamma)$, it will be seen
later that we need to compute $L_{\beta\beta}^{-1} = (I-P_{\beta\beta})^{-1}$.
The following holds for any invertible $(n-1)\times(n-1)$ matrix
$L_{\alpha,\alpha} = N^{-1} $ whose principal submatrix $L_{\beta,\beta}$ is also
invertible.
Suppose the rows/columns are ordered WLOG with $\beta=\{1,\ldots,n_1\}$,
$\gamma=\{n_1+1,\ldots,n-1\}$. Then
$${
L_{\alpha,\alpha} \cdot N 
=
\left(\apile[ll]{L_{\beta,\beta} & L_{\beta,\gamma} \\ L_{\gamma,\beta} & L_{\gamma,\gamma}}
\right)
\cdot
\left(\apile[ll]{N_{\beta,\beta} & N_{\beta,\gamma} \\ N_{\gamma,\beta} & N_{\gamma,\gamma}}
\right)
=
\left(\apile[cc]{I & 0 \\ 0 & I}\right)
}.$$
From this relation we see that
$${
\apile[llll]{
\mbox{(a)} &
L_{\beta,\beta}N_{\beta,\gamma} +  L_{\beta,\gamma} N_{\gamma,\gamma} &=& 0,
\\
\mbox{(b)} &
L_{\beta,\beta}N_{\beta,\beta} +  L_{\beta,\gamma} N_{\gamma,\beta} &=& I_{\beta,\beta}.
}
}\myeqno{equ:lll}$$
Then we can verify the formula for the inverse of $L_{\beta,\beta}$ as follows,
assuming the indicated inverses exist:
$${
\apile[lll]{
\mc{3}{l}{L_{\beta,\beta}\left[N_{\beta,\beta} -  N_{\beta,\gamma} N_{\gamma,\gamma}^{-1} N_{\gamma,\beta}\right]}
\\\mbox{\quad} &=&
I_{\beta,\beta}  -  L_{\beta,\gamma} N_{\gamma,\beta} 
- L_{\beta,\beta}  N_{\beta,\gamma} N_{\gamma,\gamma}^{-1} N_{\gamma,\beta}
\\ &=&
I_{\beta,\beta} -  L_{\beta,\gamma} N_{\gamma,\beta}
  +  L_{\beta,\gamma} N_{\gamma,\gamma} N_{\gamma,\gamma}^{-1} N_{\gamma,\beta}
\\ &=& I_{\beta,\beta}
}
}\myeqno{equ:lll2}$$
We have (where (\ref{equ:inv-schur}b) follows from (\ref{equ:lll}b)):
$${
\apile[llll]{
\mbox{(a)} &
L_{\beta,\beta}^{-1} & = & 
N_{\beta,\beta} -  N_{\beta,\gamma} N_{\gamma,\gamma}^{-1} N_{\gamma,\beta}
\\
\mbox{(b)} &
&=& N_{\beta,\beta} + L_{\beta,\beta}^{-1} L_{\beta,\gamma}
 N_{\gamma,\beta}
}
}\myeqno{equ:inv-schur}$$

\section{Applications}
{
\boldparagraph          {Hitting Time}
The expected time for a random walk starting at source \(i\) to reach target \(k\)
is
\[
H(i,k)= \sum_{j} \bm{N}(i,j,k)
= \mathrm{m}_{kk} - \mathrm{m}_{ik} + \sum_j (\mathrm{m}_{ij}-\mathrm{m}_{kj})\pi_j.
\]
The expected round-trip commute time between nodes $i$ and $k$ is
\[
\apile[r@{~}c@{~}l]{
C(i,k) &=& H(i,k) + H(k,i) =  \sum_{j} (\bm{N}(i,j,k) + \bm{N}(k,j,i))
\\
& = &
\mathrm{m}_{kk} \mathrm{m}_{ii} - \mathrm{m}_{ik} - \mathrm{m}_{ki}.
}
\]
These can easily be computed directly
from the Laplacian (see e.g. \cite{Boley11} and references therein).
\boldparagraph          {Centrality Measures}
It is easy to compute a variety of centrality measures such as
the average round-trip commute time from a given node $k$ to all other nodes
\[
\apile[r@{~~}c@{~~}l]{
C(k) &=& \sum_i C(i,k)/n = \sum_{i,j}(\bm{N}(i,j,k)+\bm{N}(k,j,i))/n
\\ &=& (\mathrm{m}_{kk} + \textsc{Trace}({\bb M}))/n ,
}
\]
where we have used $\textsc{Trace}({\bb M})=\sum_i \mathrm{m}_{ii}$ and
the fact that $\sum_i \mathrm{m}_{ki} = \sum_i \mathrm{m}_{ik} = 0$.
Another measure 
of the
importance of individual nodes in terms of bottleneck or influence
is random walk closeness [Noh and Rieger, 2004]:
\[\mbox{closeness}(k)=\sum_{i} H(i,k)= \sum_{i,j}\bm{N}(i,j,k).\]
The fundamental tensor can also be used to compute
betweenness measures such as
random walk betweenness \cite{Newman05,Kang11}:
\[\mbox{betweenness}(j) =\sum_{i\not=j,k\not=j}\Prob({i  \rightarrow   j   \rightarrow   k}),\] 
where $\Prob({i  \rightarrow   j   \rightarrow   k})$ denotes the probability of passing $j$ before reaching $k$ in a random walk starting at $i$. 
}

\boldparagraph{Probability of passage before another node}
The fundamental tensor can be used to obtain the probability of passing
through a node $j$ before reaching a node $k$, denoted
\(\Prob({i \rightarrow j  \rightarrow k})\).
This is the key concept behind the personalized hitting time trust mechanism
\cite{Buendia17}. We have the following result.

\newthm{ (Probability of ordered passage)} \label{thm:prob}
The probability of passing through $j$ on a random walk through the entire network
starting from node $i$ and before reaching $k$ is
\[ \Prob({i  \rightarrow   j   \rightarrow   k})= \bm{N}(i,j,k)\,/\,\bm{N}(j,j,k). \]

\Proof
We show WLOG 
\( \Prob({i  \rightarrow   (n-1)   \rightarrow   n}) = \bm{N}(i,n-1,n)\,/\,\bm{N}(n-1,n-1,n)\).
Consider a recurring Markov chain with transition probabilities $P$ (an $n \times n$ matrix), partitioned as
$${
P = \left[
\apile[cc|c]{Q & \bb r_1 & \bb r_2 \\
\bb s_1^T & t_{11} & t_{12} \\ \hline
\bb s_2^T & t_{21} & t_{22} }\right]
},$$
where $Q$ is $(n-2) \times (n-2)$, 
$\bb r_1, \bb r_2$ are column $(n-2)$-vectors, and
$\bb s_1^T, \bb s_2^T$ are row $(n-2)$-vectors.

\mbox{~~}

To count how many times we pass through node $n-1$ before reaching $n$,
turn node $n$ into an absorbing node and follow the prescription
in \cite{Grinstead06}.
Set $\bb s_2 = \bb 0$, $t_{21}=0$, $t_{22}=1$, and follow
\cite{Grinstead06} to find the average number of visits to any node $j$ starting from $i$
via the fundamental matrix:
$${
N = \left[
\apile[cc]{(I-Q) & -\bb r_1  \\
-\bb s_1^T & 1-t_{11}  }\right]^{-1}
= \left[\apile[cc]{W & \bb x \\ \bb y^T & z }\right],
}\myeqno{equ:countvisits}$$
partitioned conformally (so $z$ is a scalar).
In particular, \(x_i\) is the expected number of passages through node $n-1$
before reaching node $n$ when starting from node $i$,
and $z $ is the expected number of passages through node $n-1$ when starting from node
$n-1$.
So we have $x_i = \bb N(i,n-1,n)$ and $z = \bb N(n-1,n-1,n)$.

From
$${
 \left[
\apile[cc]{(I-Q) & -\bb r_1  \\
-\bb s_1^T & 1-t_{11}  }\right]
\cdot \left[\apile[cc]{W & \bb x \\ \bb y^T & z }\right]
= \left[\apile[cc]{I & \bb 0 \\ \bb 0^T & 1}\right]
}$$
we have $(I-Q)\bb x - z  \bb r_1 = \bb 0$, or
$${
\bb x = z  (I-Q)^{-1} \bb r_1 .
}\myeqno{equ:cnt}$$

Now set  $\bb s_1 = \bb 0$, $t_{11}=1$, $t_{12}=0$,
turning both $n-1,n$ into absorbing states, and follow \cite{Grinstead06}
to find the vector of probabilities of reaching $n-1$ before $n$
(starting from any node $i=1,\ldots,n-2$) as
$\bb b_1 = (I-Q)^{-1} \bb r_1 $, which from (\ref{equ:cnt}) is just
$\bb x / z  $.
\QED
\boldparagraph{Avoiding nodes}
We now consider the average hitting time to a given node
assuming we avoid a certain subset of nodes.
Specifically, we divide the set of nodes into three subsets,
$$\beta=\{1,\ldots,n_1\}, \quad \gamma=\{n_1+1,\ldots,n-1\}, \quad \delta=\{n\},$$
such that $\alpha = \{\beta,\gamma\}=\{1,\ldots,n-1\}$ in (\ref{equ:P}).
We consider the case where, starting from some node $i$ in
$\beta$, we would like to obtain the average hitting time
to node $n$ in $\delta$, {\em conditioned} on avoiding
any node in $\gamma$. We denote the expected number of passages through
an individual node $j$, starting at $i$, before reaching $n$, and
avoiding $\gamma$ as $\bm N(i,j,n,\gamma)$.
Then this desired conditional average hitting time is
$H(i,n,\gamma) = \sum_{j\in\beta}\bm N(i,j,n,\gamma)$.
Assume WLOG that the nodes are ordered so that
the probability transition
matrix is partitioned as
$${
P = 
\left(\apile[c|c]{
P_{\alpha,\alpha} & P_{\alpha,{}n} \\ \hline
P_{{}n,\alpha} &  P_{{}n{}n}
} \right)
=
\left(\apile[cc|c]{
P_{\beta,\beta} & P_{\beta,\gamma} & P_{\beta,{}n} \\
P_{\gamma,\beta} & P_{\gamma,\gamma} & P_{\gamma,{}n} \\ \hline
P_{{}n,\beta} & P_{{}n,\gamma} & P_{{}n{}n}
}\right) 
 .
}$$
If node $n$ is temporarily treated as an absorbing node, we 
obtain the $(n-1)\times(n-1)$ matrix of visit counts to any individual node in $\alpha$ before reaching
node $n$ shown in eq.\ (\ref{Nlap}).
In general,
we use the following result from \cite[chap.\ 11]{Grinstead06}:
\newthm[Proposition]\label{c0}
For $i,j\in \beta$, $\ell \in \{\gamma,n\}$:
$${
\apile[l@{~}l@{~}l]{
[(I - P_{\beta,\beta})^{-1}]_{ij} 
&=& \bm N(i,j,\{\gamma,n\})
\\ &=& \textrm{\parbox[t]{4.5in}{expected number of passages through $j$ starting from $i$ before leaving $\beta$}}
\\ {}%
[(I - P_{\beta,\beta})^{-1}\cdot P_{\beta,:}]_{i\ell} 
&=& \Prob(i \rightarrow \ell \rightarrow [\{\gamma,n\}\mbox{$-$}\{\ell\}])
\\ &=& \parbox[t]{4.5in}{\Prob(\textrm{$\ell$ is first node reached outside of $\beta$, starting from
$i$})}
}
}\label{equ:GS}$$

\newthm[Proposition]\label{c1}
For any $j\in\beta$
such that $n$ is reachable from $j$ without passing through $\gamma$, 
the average number of return visits to the 
starting node $j$
before reaching $n$ and avoiding $\gamma$ is
$${\bm N(j,j,n,\gamma) = \bm N(j,j,\{\gamma,n\}) 
=
[(I-P_{\beta,\beta})^{-1}]_{jj}
}.\myeqno{selfvisits}$$
That is, when counting return visits, the condition to avoid certain nodes makes no difference.
\Proof
Starting from $j$,
consider the three mutually exclusive events with corresponding
probabilities 
\itize{
\item
$p=\Prob(\textrm{pass $j$ before reaching $\gamma,n$})$,
\item
$q=\Prob(\textrm{reach $\gamma$ before reaching $j,n$})$,
\item
$r=\Prob(\textrm{pass $n$ before reaching $j,\gamma$})$,
}
with $p+q+r=1$.
We have the following probabilities, starting from $j$ and wandering among nodes in $\beta$ only:
\begin{itemize}
\item 
\(
a := \Prob(\textrm{pass $j$ exactly $k$ times, then reach $n$ before $\gamma$}) %
 = p^kr 
\)

\item 
\(
b := \Prob(\textrm{eventually reach $n$ before $\gamma$}) %
  = r/(q+r)=r/(1-p) = r\sum_{k=0}^{\infty}p^k
\)

\item
\(
c := \Prob(\textrm{pass $j$ exactly $k$ times, then reach $n$ or $\gamma$}) %
= p^k(1-p)
\)
\end{itemize}

Hence the probability of passing $j$ exactly $k$ times {\em conditioned} on
eventually reaching $n$ is
\begin{align*}
\Prob(\textrm{pass $j$ exactly $k$ times, given reach $n$ before $\gamma$}) &=  a/b
= c
\end{align*} 
Since the probability distributions over $k$ match, their expected values match.
\QED

\medskip

\newthm[Proposition]\label{c2}
For any $i,j\in\beta$
such that $n$ is reachable from $i$ without passing through $\gamma$,
$${{\bm N(i,j,n,\gamma) 
 =
\Prob(i \rightarrow   j \rightarrow  n \textrm{ avoiding }\gamma) 
\cdot
\bm N(j,j,\{\gamma,n\}).
}}\myeqno{condvisits}$$
The  quantity
$\Prob(i \rightarrow   j \rightarrow   n \textrm{ avoiding }\gamma)$
is the probability that a random walk starting from $i$ will visit intermediate node $j$ 
before reaching $n$, conditioned on the random walk never visiting any node in $\gamma$. 
This proposition says that the expected number of passages through $j$ is equal to the probability of reaching
$j$ initially times the expected number of return visits to $j$ once it is reached.
We have the following for the probability of reaching $j$ initially:

\newthm[Proposition]\label{c3}
For any $i,j\in\beta$
such that $n$ is reachable from $i$ without passing through $\gamma$, 
$${
{
\Prob(i \rightarrow   j \rightarrow   n \textrm{ avoiding }\gamma)
 =
\frac{\DS \Prob(i \rightarrow   j \rightarrow   \{\gamma,n\}) \cdot \Prob(j \rightarrow   n  \rightarrow \{\gamma\}) }
{\DS \Prob(i \rightarrow   n \rightarrow \{ \gamma\})}
}}.\myeqno{condprob}$$
\Proof[Proof (Sketch)]
To derive (\ref{condprob}), we temporarily treat all the nodes in $\gamma, \delta$
as absorbing states.
The numerator of the fraction is the joint probability
that starting from $i$, the random walker passes $j$ before absorption into
any of $\gamma, \delta$, and once at $j$ the walker is absorbed by
$n$ as opposed to any node in $\gamma$.
The denominator is the probability that starting from $i$,
the walker is absorbed by
$n$ as opposed to any node in $\gamma$.
\QED

We plug in the formulas from Proposition \ref{c0} 
into (\ref{condprob}) to get:
\newthm[Proposition]\label{c4}
For any $i,j\in\beta$,
$${
\apile[lllll]{
\Prob(i \rightarrow   j \rightarrow   n \textrm{ avoiding }\gamma) 
& = & 
\Prob(i \rightarrow   j \rightarrow   \{\gamma,n\})
 &\cdot& \frac{\DS \Prob(j \rightarrow   n \rightarrow \{ \gamma \} ) }
            {\DS \Prob(i \rightarrow   n \rightarrow \{ \gamma \} ) }
\\ & = &
\frac{\DS 
[(I-P_{\beta,\beta})^{-1}]_{ij}}{\DS [(I-P_{\beta,\beta})^{-1}]_{jj}}%
&\cdot&
\frac{\DS [(I-P_{\beta,\beta})^{-1}P_{\beta,:}]_{jn}%
 }
{\DS 
[(I-P_{\beta,\beta})^{-1} P_{\beta,:}]_{in}%
}
}}\myeqno{condprobformulas}$$
\QED

Propositions \ref{c2} and \ref{c3} give formulas for the reaching probabilities and average counts
conditioned on avoiding certain set of nodes.  If we plug in formulas from
Propositions \ref{c1} and \ref{c4}, we see that all the terms can be written in terms of
$(I-P_{\beta,\beta})^{-1} $ and closely related quantities.
If we wish to compute these quantities for a variety of $\beta$'s and $\gamma$'s
representing various subsets of vertices in the network, it is useful to
know how to obtain $(I-P_{\beta,\beta})^{-1} $ from the original Fundamental Tensor
$\bm N(i,j,k)$ (\ref{equ:slice}).  In the following, we assume WLOG that the nodes are ordered so
that $k = n$.

\newthm[Proposition] \label{c-scomp}
$${
\apile[l@{~}l@{~}l@{\quad \quad}l]{
\mc{4}{l}{(I-P_{\beta,\beta})^{-1}
= \bm N(\beta,\beta,\{\gamma,n\})}
\\~~ &=& 
{
[(I-P_{\alpha,\alpha})^{-1}]_{\beta,\beta}}
  & \textbf{[A]} 
\\ & & 
{
- ~
[(I-P_{\alpha,\alpha})^{-1}]_{\beta,\gamma} ~ \cdot ~
[[(I-P_{\alpha,\alpha})^{-1}]_{\gamma,\gamma}]^{-1}}
   &\textbf{[B]} %
\\ & & 
{
~ ~ ~ \cdot ~
[(I-P_{\alpha,\alpha})^{-1}]_{\gamma,\beta}}
 & \textbf{[C]} %
\\
&=& \mc{2}{l}{
\underbrace{\bm N(\beta,\beta,n)}_{\textbf{[A]}}%
-
\underbrace{\bm N(\beta,\gamma,n)%
\bm N(\gamma,\gamma,n)^{-1}}_{\textbf{[B]}}%
\underbrace{\bm N(\gamma,\beta,n)}_{\textbf{[C]}}%
}
\\
&=&
[(I-P_{\alpha,\alpha})^{-1}]_{\beta,\beta}
  & \textbf{[A]} %
\\ &&
- ~
{(I-P_{\beta,\beta})^{-1} (P_{\beta,\gamma})}
   &\textbf{[B]} %
\\ && \mbox{~~} \cdot ~
{[(I-P_{\alpha,\alpha})^{-1}]_{\gamma,\beta}}
 & \textbf{[C]} %
}
}$$
where the individual entries in the last expression represent
$${
\apile[l@{~}l@{~}l]{
\textbf{[A]}_{ij} &=& \textrm{\parbox[t]{4.0in}{expected number of passages through node $j\in\beta$ starting at $i\in\beta$}} \\
\textbf{[B]}_{i\ell} &=& \parbox[t]{4.0in}{\Prob(\textrm{$\ell${$\in$}$\gamma$ is the first node visited outside of $\beta$, starting at $i${$\in$}$\beta$)}}
\\ {}%
\textbf{[C]}_{\ell j} &=& \textrm{\parbox[t]{4.0in}{expected number of passages through node $j\in\beta$ starting at $l\in\gamma$}} 
}
}$$
In all three cases the walks stop when they reach $n$.
Here the notation $\bm N(\beta,\gamma,n)$ denotes the sub-matrix of
$\bm N(:,:,n)$ consisting of the rows indexed by $\beta$
and columns indexed by $\gamma$.
\Proof[Proof (Sketch)]
The first equality is just a rewrite of (\ref{equ:inv-schur}a) and the last
equality is a rewrite of (\ref{equ:inv-schur}b).
The interpretation of $\textbf{[B]}_{i\ell}$ as the indicated probability comes from
Proposition \ \ref{c0} and (\ref{equ:inv-schur}b).
\QED

In the above formulas, if the avoidance set $\gamma$ is relatively small, then
the cost of the above formulas will be modest, at most quadratic in $n$ once
$(I-P_{\alpha,\alpha})^{-1}$ has been obtained.

\section{Examples}
\boldparagraph{Illustrative Example}
We give a small example with four nodes.
{\def\c#1{\zz[c]{\LARGE$\bigcirc$}$\!#1$}
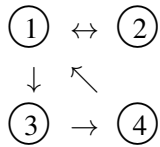
$${
\apile[c@{\quad}c@{\quad\,\,}c]{
\c1 & \leftrightarrow & \c2 \\[1ex]
\downarrow  & \nwarrow & \\[1ex]
\c3 & \rightarrow & \c4
}
}$$
}%
We define a visit to a node $j$ as one departure from that node.
Hence the visit counts for target nodes are zero:
\(\bm{N}(i,j,k) = 0\) for \(i = k\) or \(j=k \).
In this case the probability transition matrix and normalized Laplacian are
$${
P = \begin{pmatrix}
0      & \frc12  & \frc12  & 0     \\
1      & 0       & 0       & 0     \\
0      & 0       & 0       & 1     \\
1      & 0       & 0       & 0     
\end{pmatrix},
\quad
L = \begin{pmatrix}
1      & -\frc12 & -\frc12 & 0     \\
-1     & 1       & 0       & 0     \\
0      & 0       & 1       & -1    \\
-1     & 0       & 0       & 1     
\end{pmatrix},
}$$
$${
L^+ = \frac{1}{28} \left(\apile[rrrr]{
8      &  -3     & -3      & -10   \\
0      & 21      & -7      & -14   \\
-8     & -11     & 17      & 10    \\
0      & -7      & -7      & 14    
}\right),
}$$ 
where the recurring probabilities are $\bm{\pi}$
$=$
$( 0.4 ;\, 0.2 ;\,0.2 ;\, 0.2 )$.
The computed tensor is then
$${ 
\apile[ll]{
\bm{N}_{::1} = 
\begin{pmatrix}
0 & 0 & 0 & 0 \\
0 & 1 & 0 & 0 \\
0 & 0 & 1 & 1 \\
0 & 0 & 0 & 1 
\end{pmatrix} 
&
\bm{N}_{::2} = 
\begin{pmatrix}
2 & 0 & 1 & 1 \\
0 & 0 & 0 & 0 \\
2 & 0 & 2 & 2 \\
2 & 0 & 1 & 2 
\end{pmatrix}
}
}$$
$${
\apile[ll]{
\bm{N}_{::3} = 
\begin{pmatrix}
2 & 1 & 0 & 0 \\
2 & 2 & 0 & 0 \\
0 & 0 & 0 & 0 \\
2 & 1 & 0 & 1 
\end{pmatrix} 
&
\bm{N}_{::4} = 
\begin{pmatrix}
2 & 1 & 1 & 0 \\
2 & 2 & 1 & 0 \\
0 & 0 & 1 & 0 \\
0 & 0 & 0 & 0
\end{pmatrix}.
}
}$$

\begin{figure}[t]
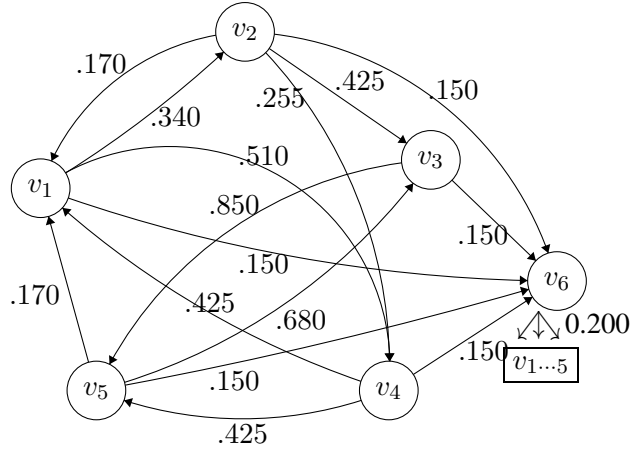

\centering%
\usebox{\GGaug}
\caption{Sample network for trust mechanism model. Node $v_6$ is
used to model the global evaporation node.}
\label{fig:trust}
\end{figure}

\boldparagraph{Trust Mechanism}
We consider the personalized hitting time (PHT) trust mechanism, which is one way to model the propagation of trust through
a network of actors \cite{Liu16,Buendia17}.
Consider the network of Fig.\ \ref{fig:trust} with probability transition matrix
{\def\1{\mc{1}{c}{0}}
$${
P = \left(\apile[r@{~~~}r@{~~~}r@{~~~}r@{~~~}r@{~~~}r]{
       \1  &  0.340  &     \1  &  0.510  &     \1  &  0.150  \\
    0.170  &     \1  &  0.425  &  0.255  &     \1  &  0.150  \\
       \1  &     \1  &     \1  &     \1  &  0.850  &  0.150  \\
    0.425  &     \1  &     \1  &     \1  &  0.425  &  0.150  \\
    0.170  &     \1  &  0.680  &     \1  &     \1  &  0.150  \\
    0.200  &  0.200  &  0.200  &  0.200  &  0.200  &     \1
}\right)%
}\myeqno{eq:trust}$$
}%
Here we seek to infer trust values for every node from the viewpoint of every other node,
where $v_6$ is an outside evaporation node to enforce proximity in the trust values
(See \cite{Liu16,Buendia17} for details).  The personalized hitting time trust mechanism
value is defined as
$${\PHT(i,j)= \Prob(i \rightarrow j \rightarrow 6 ).}$$
According to Theorem \ref{thm:prob}, this can be computed as
${
\PHT(i,j) = \bm{N}(i,j,6)\,/\,\bm{N}(j,j,6) .
}$
For example, from the viewpoint of node 4, the trust values for nodes $1,\ldots,5$ are
$${
\PHT(4,1:5) = ( 0.5962,\, 0.2913,\, 0.5332 ,\, 1.0 ,\, 0.6573),
}$$
where node 5 enjoys the highest trust value from the point of view of node 4.
If we instead consider only paths that avoid node 2 (for instance, node 2 is known to be
a bad actor), we can compute the PHT values based on the restricted matrix of counts
$\bm{N}(:,:,6,2)$ instead of $\bm{N}(:,:,6)$, to obtain the modified trust values
$${
\PHT_2(4,1:5) = ( 0.5962 ,\, 0.0000 ,\,0.3872 ,\,1.0 ,\,0.5426 ),
}$$
where it is seen that the node with the highest trust value has changed to node 1.


\end{document}